\documentclass[11pt]{article}
\usepackage{graphicx}
\usepackage{latexsym}
\usepackage{amsmath}
\usepackage{amssymb}
\usepackage{hyperref}
\usepackage{amsfonts}
\usepackage{color,xcolor}
\usepackage{manfnt}
\usepackage{float}

\newcommand{\qed}{\nobreak \ifvmode \relax \else
      \ifdim\lastskip<1.5em \hskip-\lastskip
      \hskip1.5em plus0em minus0.5em \fi \nobreak
      \vrule height0.75em width0.5em depth0.25em\fi}

\def \Vh0{\stackrel{\circ}{V}_h} 
   
\def\Om{\Omega}   
\newcommand{\q}{\quad}    
  \def\f{\frac}

  \def\x{{\bf x}}

\def\bb{\begin{equation}} \def\ee{\end{equation}}
\def\beqn{\begin{eqnarray}}  \def\eqn{\end{eqnarray}}
\def\beqnx{\begin{eqnarray*}} \def\eqnx{\end{eqnarray*}}

\title{Two-level space-time domain decomposition methods for unsteady inverse problems}
\author{Xiaomao Deng \footnote{
Laboratory for Engineering and Scientific Computing,
        Shenzhen Institutes of Advanced Technology, Chinese Academy of Sciences,
        Shenzhen, Guangdong 518055, {P. R. China}.
        (xm.deng@siat.ac.cn)}
\quad\quad Xiao-chuan Cai\footnote{Department of Computer Science, University of Colorado Boulder, Boulder, CO 80309, {USA}.
(cai@colorado.edu)}
\quad\quad Jun Zou\footnote{Department of Mathematics, The Chinese University of Hong Kong, Shatin N.T., Hong Kong, {P. R. China}. The work of this author was substantially supported by
Hong Kong RGC grants (projects 14306814 and 405513). (zou@math.cuhk.edu.hk)}}

\date{}
\begin{document}
\maketitle
\begin{abstract}
As the number of processor cores on supercomputers becomes larger and larger, algorithms with high degree of parallelism attract more attention.
In this work, we propose a novel space-time coupled algorithm for solving an inverse problem associated with the
time-dependent convection-diffusion equation in three dimensions.
We introduce a mixed finite element/finite difference method and a one-level and a two-level space-time parallel domain decomposition preconditioner for the Karush-Kuhn-Tucker (KKT) system induced from reformulating the inverse problem as an output least-squares optimization problem in the space-time domain.
The new full space approach eliminates the sequential steps of the optimization outer loop and the inner forward and backward time marching processes, thus achieves high degree of parallelism.
Numerical experiments validate that this approach is effective and robust for recovering unsteady moving sources.
We report strong scalability results obtained on a supercomputer with more than 1,000 processors.

{\bf Keywords: }Space-time method;~Domain decomposition method;~Unsteady inverse problem;~Pollutant source identification;~Parallel computing
\end{abstract}
\section{Introduction}
In this paper, we consider an inverse problem associated with the time-dependent convection-diffusion equation
defined in $\Om\in \mathbf{R}^3$:
\begin{equation}
\begin{cases}
\displaystyle\frac{\partial C}{\partial t}=\nabla\cdot(a(\mathbf{x})\nabla C)-\nabla\cdot (\mathbf{v}(\mathbf{x})C)+ f(\x, t), \quad 0< t< T\,, ~\mathbf{x}\in \Omega \\
C(\mathbf{x},t)=p(\mathbf{x},t),\quad \x\in \Gamma_1 \q\\
a(\mathbf{x})\displaystyle\frac{\partial C}{\partial \mathbf{n}}=q(\mathbf{x},t),\quad \x\in \Gamma_2 \q\\
C(\mathbf{x},0)=C_0(\mathbf{x}), \q \x\in \Om\,,
\label{eq:cfeq}
\end{cases}
\end{equation}
where $f(\x, t)$ is the source term to be recovered, $a(\mathbf{x})$ and
$\mathbf{v}(\mathbf{x})$ are the given diffusivity and convective coefficients.
$\Gamma_1$ and $\Gamma_2$ are two disjoint parts of the boundary $\partial\Om$.
Dirichlet and Neumann boundary conditions are imposed respectively
on $\Gamma_1$ and $\Gamma_2$.
When the observation data $C(x,t)$ is available at certain locations,
several classes of inverse problems associated with the convection-diffusion equation (\ref{eq:cfeq})
have been investigated, such as the recovery of the diffusivity coefficient with applications in, for examples, laminar wavy film flows \cite{KGM11}, and  flows in porous media  \cite{NKM09}, the recovery of the source with applications in, for examples, convective heat transfer problems \cite{MD09}, indoor airborne  pollutant tracking \cite{LZ07}, ground water contamination modeling \cite{RR05,SK94,SK95,YSL08}, etc.

The main focus of this work is to study the following inverse problem:
given the measurement data $C^{\epsilon}(\mathbf{x},t)$ of $C(\mathbf{x},t)$
at some locations inside $\Omega$ for the period $0<t<T$ ($\epsilon$ denotes the noise level),
we try to recover the time-varying source locations and intensities, i.e., the source function $f(\x, t)$ in equation (\ref{eq:cfeq}).
In the last decades,
several types of numerical methods have been developed for retracing the sources, such as the
explicit method \cite{H09}, the quasi-reversibility method \cite{SK95}, the statistical method \cite{SK97} and the Tikhonov optimization method \cite{AB01,DZZ13,GER83,WY10}.
Among these methods, the Tikhonov optimization method is the most popular one,
which reformulates the original inverse source problem into an output least-squares optimization problem with PDE-constraints,
and by including appropriate regularizations
it ensures the stability of the resulting optimization problem \cite{EHN98,W03}.
Various techniques are available for solving the induced first-order optimality system \cite{AB03,AB01,B96,BG99,KT51}.

We define the following objective functional with Tikhonov optimization:
\begin{equation}
J(f)=\displaystyle\frac{1}{2}\int_{0}^T \int_{\Omega} A(\x)(C(\mathbf{x}, t)-C^{\epsilon}(\mathbf{x},t))^2\,d\x dt  + N_{\beta}(f)\,,
\label{eq:optim}
\end{equation}
where $A(\x)$ is the data range indicator function, namely $A(\x)=\sum_{i=1}^{s} \delta(\x-\x_i)$,
and $\x_1$, $\x_2$, $\cdots$, $\x_{s}$ are a set of specified locations, where
the concentration $C(\mathbf{x}, t)$ is measured and denoted by $C^{\epsilon}(\mathbf{x_i},t)$.
The term $N_{\beta}(f)$ in (\ref{eq:optim}) is called the regularization with respect to the source.
Since $f(\x, t)$ depends on both space and time,
we propose the following space-time $H^1$-$H^1$ regularization:
\begin{equation}
 N_{\beta}(f)=\displaystyle\frac{\beta_1}{2}\int_0^T\int_{\Om} | \dot{f}(\x, t)|^2 d\x dt + \displaystyle\frac{\beta_2}{2}\int_0^T\int_{\Om} |\nabla f|^2 d\x  dt\,.\label{eq:reguterm}
\end{equation}
Here $\beta_1$ and $\beta_2$ are two regularization parameters. Other regularizations, such as $H^1$-$L^2$,
may be used, but we will show later by numerical experiments that $H^1$-$H^1$
regularization may offer better numerical reconstructions.

Traditionally the problem (\ref{eq:cfeq})-(\ref{eq:reguterm}) is solved by a reduced space sequential quadratic programming(SQP) method \cite{DZZ13,WY10},
which can be described as follows:
\begin{enumerate}
\item[] {\bf Optimization loop} (sequential)
\item[] Step 1: Solve a forward-in-time state equation
   \begin{enumerate}
   \item[] Loop in time (sequential)\\
   - Solve the steady-state equation for this time step (parallel)
   \item[] End loop
   \end{enumerate}
\item[] Step 2: Solve a backward-in-time adjoint equation
   \begin{enumerate}
   \item[] Loop in time (sequential)\\
   - Solve the steady-state equation for this time step (parallel)
   \item[] End loop
   \end{enumerate}
\item[] Step 3: Solve objective equation (parallel)
\item[] End loop
\end{enumerate}

Parallelization strategies for solving the problem with reduced space SQP methods include
by keeping the sequential steps of the outer loop and applying parallel-in-space algorithms such as domain decomposition, multigrid methods to the subsystems at each time step \cite{AB05}.
Reduced space SQP methods for unsteady inverse problems needs repeatedly solving the state equation, the adjoint equation and the objective equation, thus it divides the problem into many subproblems, the memory cost is low.
However reduced space SQP methods sometimes are quite time-consuming to achieve convergence.
Because of the sequential steps in the optimization loop and in the forward and backward time-marching processes, it is less ideal for parallel computers with a large number of processor cores compared to full space SQP methods.
Full space method in \cite{CLZ09,CC12} has been studied for steady state problems, for unsteady problems,
it needs to eliminate the sequential steps in the outer loop and solve the full space-time problem as a coupled system.
Because of the much larger size of the system, the full space approach may not be suitable for small computer systems, but it has fewer sequential steps and thus offers a much higher degree of parallelism required by large scale supercomputers.

Finding suitable parallelization strategies for the optimization loop and the inner time loops is an active research area.
An unsteady PDE-constrained optimization problem was solved in \cite{YPC12} for the boundary control of unsteady incompressible flows by solving a subproblem at each time step.
It has the sequential time-marching process and each subproblem is steady-state.
The parareal algorithms were studied in \cite{BBMTZ02,FC03,LMT01},
which involve a coarse (coarse mesh in the time dimension) solver for prediction and a fine (fine mesh in the time dimension) solver for correction.
Parallel implicit time integrator method (PITA), space-time multigrid, multiple shooting methods can be categorized as improved versions of the parareal algorithm \cite{GP08,GV07}. The parareal algorithm combined with domain decomposition method \cite{MT05} or multigrid method can be powerful. So far, most references on parareal related studies focus mainly on the stability and convergence \cite{GH08}.

In this paper, we propose a fully implicit, mixed finite element and finite difference discretization scheme for the continuous KKT system, and a corresponding space-time overlapping Schwarz preconditioned solver for the unsteady inverse source identification problem in three dimensions.
The method removes all the sequential inner time steps and achieves full parallelization in both space and time.
We study the most general form of the source function, in other words, we reconstruct the time history of the distribution and intensity profile of the source simultaneously.
Furthermore, to resolve the dilemma that the number of linear iterations of one-level methods increases with the number of processors, we develop a two-level space-time hybrid Schwarz preconditioner which offers better  performance in terms of the number of iterations and the total compute time.

The rest of the paper is arranged as follows. In Section~\ref{sec:Strong formulation of KKT system}, we describe
the mathematical formulation of the inverse problem and the derivation of the KKT system.
We propose, in Section~\ref{sec:scalable solvers}, the main algorithm of the paper, and discuss several technical issues involved in the fully implicit discretization scheme and the one- and two-level overlapping Schwarz methods for solving the KKT system.
Numerical experiments for the recovery of 3D sources are given in Section~\ref{sec:Numerical Examples},
 and  some concluding remarks are provided in Section~\ref{sec:conclusions}.
\section{Strong formulation of KKT system}\label{sec:Strong formulation of KKT system}
We formally write (\ref{eq:cfeq}) as an operator equation
$L(C, f)=0$.
For $G\in H^{1}(\Om)$, the following Lagrange functional \cite{AB03,KZ98} transforms the PDE-constrained optimization problem (\ref{eq:optim}) into an unconstrained minimization problem. Let
\begin{equation}
\begin{split}
\mathcal{J}(C,f,G)&=\displaystyle\frac{1}{2} \int_{0}^T \int_{\Omega} A(\x)(C(\x,t)-C^{\epsilon}(\x, t))^2d\x d t\\
&+ N_{\beta}(f)+ (G,L(C,f))\,,
\end{split}
\label{eq:lagrange1}
\end{equation}
where $G$ is a Lagrange multiplier or adjoint variable, and $(G,L(C,f))$ denotes their inner product.
Two approaches are available for solving (\ref{eq:lagrange1}), the optimize-then-discretize approach and the discretize-then-optimize approach.
The first approach derives a continuous optimality condition system and then applies certain discretization scheme, such as a finite element method to obtain a discrete system ready for computation.
The second approach discretizes the optimization function $\mathcal{J}$, and then the objective functional becomes a finite dimensional quadratic polynomial.  The solution algorithm is then based on the polynomial system.
The two approaches perform the approximation and discretization at different stages, both have been applied successfully  \cite{PBC06}, we use the optimize-then-discretize approach in this paper.

The first-order optimality conditions for (\ref{eq:lagrange1}), i.e., the KKT system,
is obtained by taking the variations with respect to $G$, $C$ and $f$ as
\begin{equation}
\begin{cases}
\mathcal{J}_G(C,f,G)v = 0\\
\mathcal{J}_C(C,f, G)w = 0\\
\mathcal{J}_{f}(C,f, G)g = 0
\end{cases}
\label{eq:kkt}
\end{equation}
for all $v, w\in L^2(0,T; H^1_{\Gamma_1}(\Om))$ with zero traces on $\Gamma_1$ and $g\in   H^1(0,T; H^1(\Om))$.

Using integration by part, we obtain the strong form of the KKT system:
\begin{equation}
\begin{cases}
\displaystyle\frac{\partial C}{\partial t}-\nabla \cdot(a\nabla C) +\nabla\cdot (\textbf{v}(\x) C)-f =0 \\
-\displaystyle\frac{\partial G}{\partial t}-\nabla \cdot(a\nabla G) -\textbf{v}(\x)\cdot\nabla G +A(\x)C = A(\x)C^{\epsilon}\\
G+\beta_1\displaystyle\frac{\partial^2 f}{\partial t^2}+\beta_2 \Delta f =0\,.
\end{cases}
\label{eq:kkt2}
\end{equation}
To derive the boundary, initial and terminal conditions for each variable of the equations, we make use of the property that (\ref{eq:kkt}) holds for arbitrary directional functions $v, w$ and $g$. For the state equation (i.e. the first one of (\ref{eq:kkt}) or (\ref{eq:kkt2})) it is obvious to maintain the same conditions given by  (\ref{eq:cfeq}). For the adjoint equation (the second one of (\ref{eq:kkt}) or (\ref{eq:kkt2})), by multiplying the test function $w\in L^2(0,T; H^1_{\Gamma_1}(\Om))$ with $w(\cdot, 0)=0$, $a(\x)\frac{\partial w}{\partial \mathbf{n}}=0$ on $\Gamma_2$,  we have
\begin{align*}
\mathcal{J}_{C}(C,f, G)w &=\int_{0}^T\int_{\Omega}A(\x)(C(\mathbf{x},t)-C^{\epsilon}(\mathbf{x},t)) wd\x dt+\int_{\Omega} G(\x,T)w(\x,T) d \x\\
& -\int_0^T\int_{\Omega} \left(\displaystyle\frac{\partial G}{\partial t}+\nabla\cdot (a(\mathbf{x})\nabla G)+ \mathbf{v}(\x)\cdot \nabla G\right)w d \x d t\\
&-\int_0^T\int_{\Gamma_1} \left(a(\mathbf{x})\displaystyle\frac{\partial w}{\partial \mathbf{n}}\right)G d \Gamma d t\\
& +\int_0^T\int_{\Gamma_2} \left(a(\mathbf{x})\displaystyle\frac{\partial G}{\partial \mathbf{n}} + \mathbf{v}(\x)\cdot \mathbf{n} \right)w d \Gamma d t\,.
\end{align*}
By the arbitrariness of $w$, the boundary and terminal conditions for $G$ are derived:
\begin{align*}
&G(\x,t)=0, \quad \x\in \Gamma_1,~t\in [0,T] \\
&a(\mathbf{x})\displaystyle\frac{\partial G}{\partial \mathbf{n}} + \mathbf{v}(\x)\cdot \mathbf{n} = 0, \quad \x\in \Gamma_2,~t\in [0,T] \\
&G(\x,T) =0, \quad \x\in \Om\,.
\end{align*}
Similarly for the third equation of (\ref{eq:kkt}) or (\ref{eq:kkt2}), we can deduce
\begin{align*}
\mathcal{J}_{f}(C,f, G)g &=-\int_{0}^T\int_{\Omega}G g d\x dt +\int_0^T\int_{\Omega} (\dot{f} \dot{g}+\nabla f\cdot \nabla g) d \x d t\\
&=-\int_{0}^T\int_{\Omega}G g d\x dt +(\dot{f} g)|_{t=0, T}-\int_{0}^T\int_{\Omega}\ddot{f}g \\
&\mbox{\quad~} + \int_0^T\int_{\partial\Omega}
\displaystyle\frac{\partial f}{\partial \mathbf{n}}g d \Gamma d t - \int_0^T\int_{\Omega} \Delta f g d \x d t \\
&=-\int_{0}^T\int_{\Omega} (G +\ddot{f}+\Delta f)g d \x d t + (\dot{f} g)|_{t=0, T} + \int_0^T\int_{\partial\Omega} \displaystyle\frac{\partial f}{\partial \mathbf{n}}g d \Gamma d t\,.
\end{align*}
Using the arbitrariness of $g$, we derive the boundary, initial and terminal conditions for $f$:
\begin{equation}
\displaystyle\frac{\partial f}{\partial t} =0 \quad  \mbox{for}  ~~t=0, T,~\x\in \Om\,;
\quad \quad \displaystyle\frac{\partial f}{\partial \mathbf{n}} =0 \quad \mbox{for} ~~\x\in \partial\Om,~t\in [0, T]\,.
\label{eq:bdyf}
\end{equation}

\section{A fully implicit and fully coupled method}\label{sec:scalable solvers}
In this section, we first introduce a mixed finite element and finite difference method for the discretization of the continuous KKT system derived in the previous section, then we briefly mention the algebraic structure of the discrete system of equations. In the second part of the section, we introduce the one- and two-level space-time Schwarz preconditioners that are the most important components for the success of the overall algorithm.
\subsection{Fully-implicit space-time discretization}\label{sec:Discretization Scheme}
In this subsection, we introduce a fully-implicit finite element/finite difference scheme to discretize (\ref{eq:kkt2}).
To discretize the state and the adjoint equations, we use a second-order Crank-Nicolson finite difference scheme in time and a piece-wise linear continuous finite element method in space.
Consider a regular triangulation $\mathcal{T}^h$ of domain $\Omega$, and a time partition
$P^{\tau}$ of the interval $[0, T]$: $0=t^0<t^1<\cdots< t^M=T,$ with $t^n=n \tau, \tau=T/M$.
Let $V^h$ be the  piecewise linear continuous finite element space on $\mathcal{T}^h$,
and  $\mathring{V}^h$ be the subspace of $V^h$ with zero trace on $\Gamma_1$.
We introduce the difference quotient and the averaging of a function $\psi(\x, t)$ as
\[\partial_{\tau} \psi^n (\x)= \displaystyle\frac{\psi^n(\x)-\psi^{n-1}(\x)}{\tau}, \quad \bar{\psi}^n(\x)=\displaystyle\frac{1}{\tau}\int_{t^{n-1}}^{t^n} \psi(\x,t) dt\,,\]
with $\psi^n(\x):=\psi(\x, t^n)$.
Let $\pi_h$ be the finite element interpolation associated with the space $V^h$,
then we obtain the  discretizations for the state and adjoint equations by finding the sequence of approximations $C_h^n, G_h^n\in V^h$, such that $C_h^0=\pi_hC_0$, $G_h^M=\mathbf{0}$,  and $C_h^n(\x)=\pi_h p(\x, t^n), G_h^n(\x)=0$ for $\x\in \Gamma_1$, and satisfying
\begin{equation}
\begin{cases}\label{eq:kkt3}
(\partial_{\tau} C_h^n, v_h) + (a\nabla  \bar{C}_h^n, \nabla v_h) + (\nabla\cdot (\mathbf{v}\bar{C}_h^n),v_h) =  (\bar{f}^n_h, v_h)+\langle\bar{q}^n,v_h\rangle_{\Gamma_2}, ~~\forall\,v_h\in \mathring{V}^h\,\\
-(\partial_{\tau} G_h^n, w_h) + (a\nabla  \bar{G}_h^n, \nabla w_h) + (\nabla\cdot (\mathbf{v}w_h),\bar{G}_h^n)
\\
=-(A(\x)(\bar{C}_h^{n}(\x,t)-\bar{C}^{\epsilon,n}(\x, t)), w_h), ~~\forall w_h\in \mathring{V}^h\,.\\
\end{cases}
\end{equation}

Unlike the approximations of  the forward and adjoint equations in (\ref{eq:kkt3}),
we shall approximate the source function $f$ differently.
We know that the source function satisfies
an elliptic equation (see the third equation in (\ref{eq:kkt2})) in the space-time domain $\Om \times (0,T)$.
So we shall apply $\mathcal{T}^h \times P^{\tau}$ to generate a partition of the space-time domain
$\Om \times (0,T)$, and then  apply
the piecewise linear finite element method in both space (three dimensions) and time (one dimension),
denoted by $W_h^{\tau}$,
to approximate the source function $f$.
Then the equation for $f\in W_h^{\tau}$ can be discretized as follows:
Find the sequence of $f_h^n$ for $n=0, 1$, $\cdots$, $M$ such that
\begin{equation} \label{eq:kkt6}
-  (G_h^n, g_h^{\tau}) +\beta_1( \partial_{\tau} f_h^n,\partial_{\tau} g_h^{\tau}) +  \beta_2 (\nabla f_h^n, \nabla  g_h^{\tau})=0,  ~~\forall\,g_h^{\tau}\in W_h^{\tau}\,.
\end{equation}
The coupled system (\ref{eq:kkt3})-(\ref{eq:kkt6}) is the so-called fully discretized  KKT system.
In the Appendix, we provide some details of the  discrete structure of this KKT system.

\subsection{One- and two-level space-time Schwarz preconditioning}\label{sec:space-time schwarz}
Usually, the unknowns of the KKT system (\ref{eq:kkt3})-(\ref{eq:kkt6}) are ordered physical variable by physical variable,
namely in the form
$$
\tilde U=(C^0, C^1, \cdots, C^{M}, G^0, G^1, \cdots, G^{M}, f^0, f^1, \cdots, f^{M})^T\,.
$$
Such ordering are used extensively in SQP methods \cite{DZZ13}.
In our all-at-once method, the unknowns $C,G$ and $f$ are ordered mesh point by
mesh point and time step by time step, and all unknowns associated with a point stay together as a block.
At each mesh point $\x_{j}$, $j=1$, $\cdots$, $N$, and time step $t^n$, $n=0$, $\cdots$, $M$, the unknowns are arranged in the order of $C_{j}^n, G_{j}^n,f_{j}^n$.
Such ordering avoids zero values on the main diagonal of the matrix and has better cache performance for point-block LU (or ILU) factorization based subdomain solvers. More precisely, we define the solution vector
\begin{eqnarray*}
U&=&(C_{1}^0,G_{1}^0, f_1^0,\cdots, C_{N}^0, G_{N}^0,f_N^0, C_{1}^1, G_{1}^1,f_{1}^1,\cdots,C_{N}^1, G_{N}^1,f_{N}^1,\cdots, C_{1}^M,G_{1}^M,f_{1}^M, \\
 &&\cdots,C_{N}^M, G_{N}^M,f_N^M)^T\,.
\end{eqnarray*}
then the linear system (\ref{eq:kkt3})-(\ref{eq:kkt6})  is rewritten as
\begin{equation}\label{eq:kktsys2}
{\,\ } \hskip-1truecm F_h U=b\,,
\end{equation}
where $F_h$ is a sparse block matrix of size $(M+1)(3N)$ by $(M+1)(3N)$ with the following block structure:
\[{\,\ } \hskip-1truecm F_h=\left(\begin{array}{ccccc}S_{00}&S_{01}&\mathbf{0}&\cdots&\mathbf{0}\\ S_{10}&S_{11}&S_{12}&\cdots&\mathbf{0}\\ \mathbf{0} & \ddots & \ddots  &\ddots &\mathbf{0} \\
\mathbf{0}&\cdots &S_{M-1,M-2}&S_{{\tiny M-1,M-1}}&S_{{\tiny M-1,M}}\\ \mathbf{0}&\cdots&\mathbf{0}&S_{M,M-1}&S_{M,M}\end{array}\right)\,,\]
where the block matrices $S_{ij}$ for $0\leq i,j \leq M$ are of size $3N\times 3N$ and most of its elements are zero matrices except the ones in the tridiagonal stripes $\{S_{i,i-1}\}, \{S_{i,i}\}, \{S_{i,i+1}\}$.
It is noted that if we denote the submatrices for $C, G$ and $f$ of size $N\times N$ respectively by $S_{ij}^C, S_{ij}^G
, S_{ij}^f$ in each block $S_{ij}$, the sparsity  of the matrices are inconsistent, namely,
 $S_{ij}^f$ is the densest and $S_{ij}^G$ is the sparest.
 This is due to the discretization scheme we have used.
The system (\ref{eq:kktsys2}) is large-scale and ill-conditioned, therefore is difficult to solve because the space-time coupled system is denser than the decoupled system, especially for three dimensional problems.
We shall design the preconditioner by extending the classical spatial Schwarz preconditioner to include both spaial and temporal variables. Such an approach eliminates all sequential steps and the unknowns at all time steps are solved simultaneously.
We use a right-preconditioned Krylov subspace method to solve (\ref{eq:kktsys2}),
 \[{\,\ } \hskip1truecm F_h M^{-1}U'=b\,,\]
where $M^{-1}$ is a space-time Schwarz preconditioner and $U=M^{-1}U'$.

\begin{figure}
\centering{
\includegraphics[width=6cm]{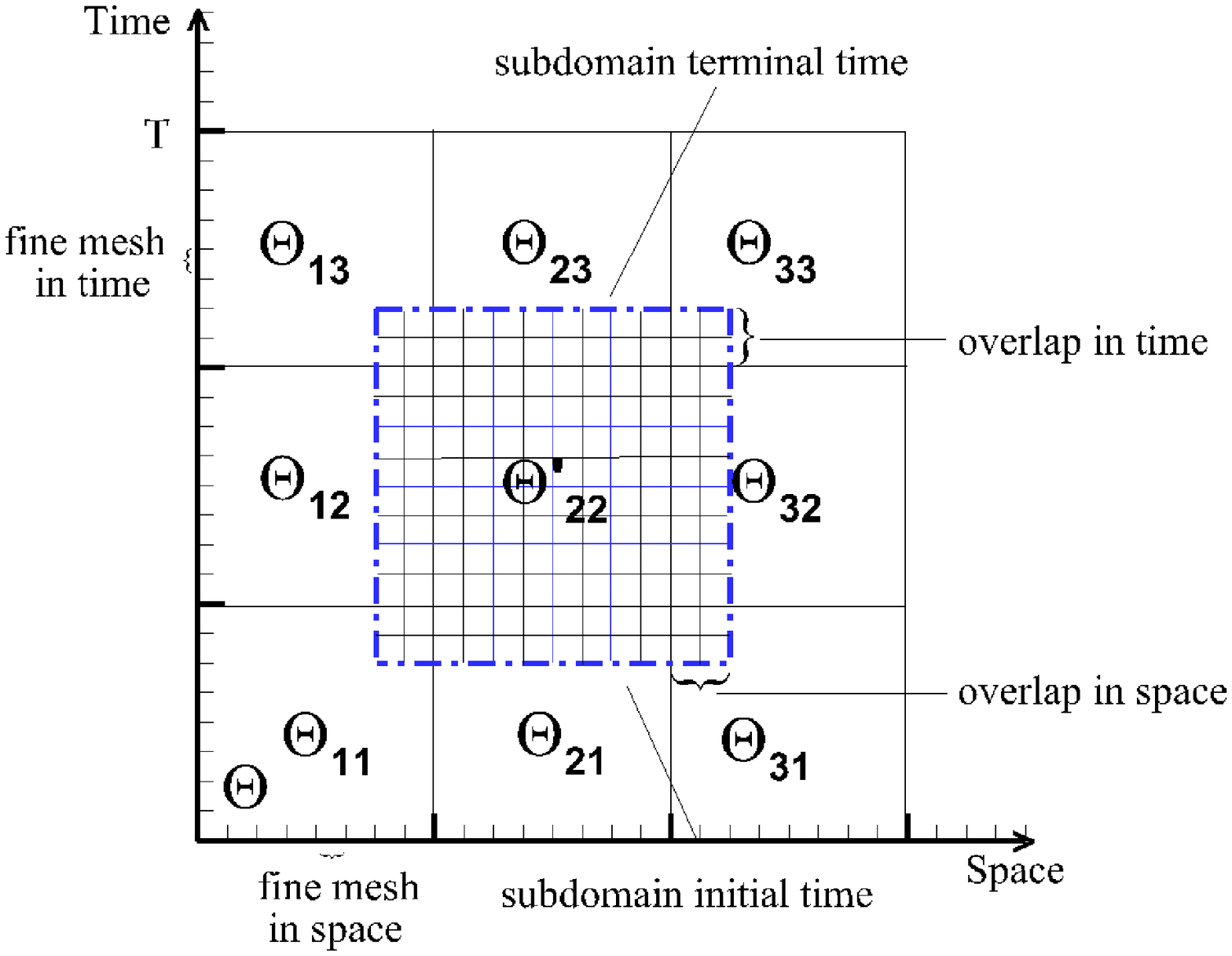}
\includegraphics[width=6cm]{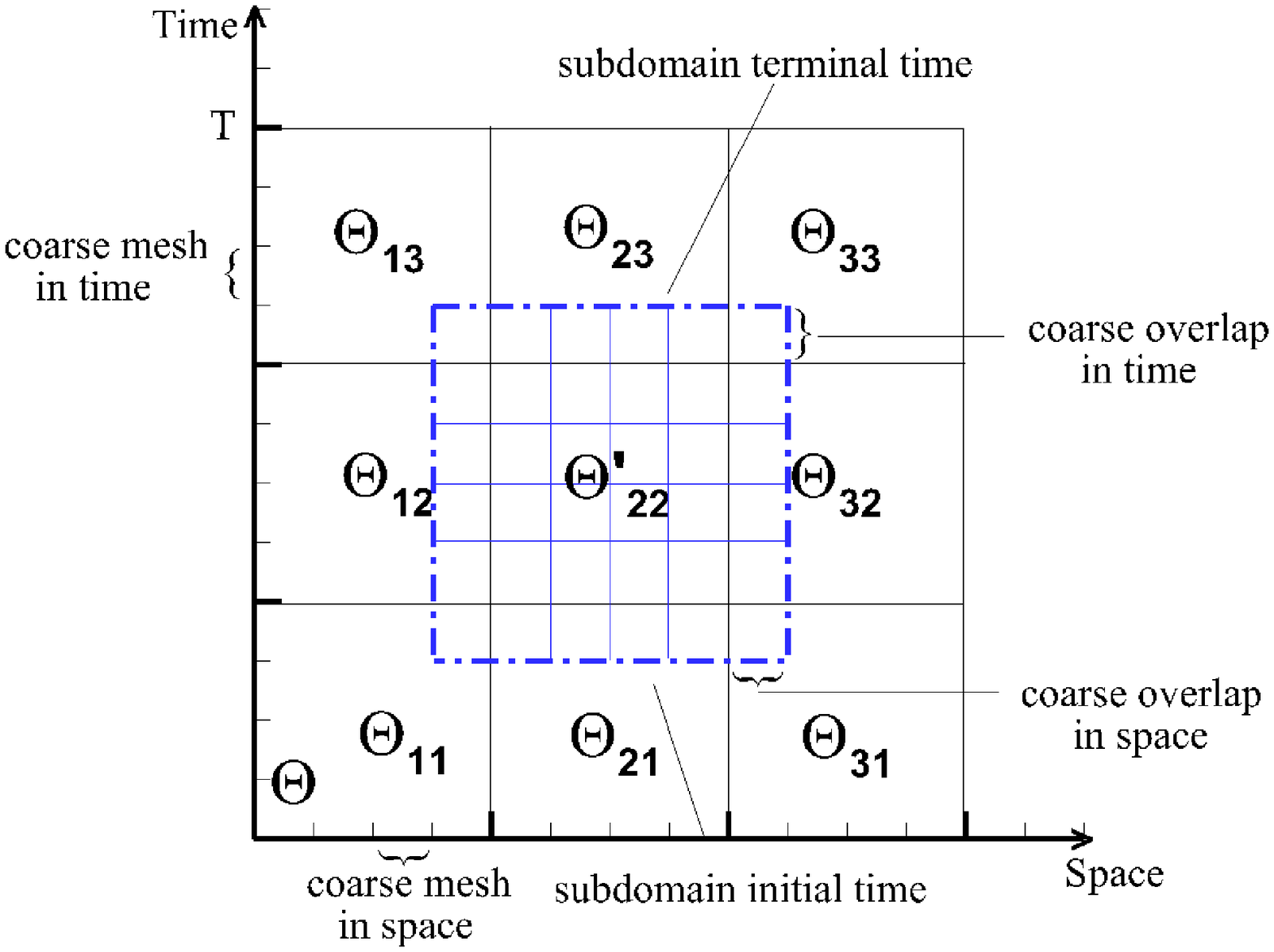}
\vspace{-0.2in}
\caption{Left: a sample overlapping decomposition in space-time domain $\Theta$ on a fine mesh. Right: the same decomposition on the coarse mesh. }\label{fig:ddm}
}
\end{figure}

Denoting the space-time domain by $\Theta=\Om\times (0,T)$, an overlapping decomposition of $\Theta$ is defined as follows:
we divide $\Om$ into $N_s$ subdomains, $\Om_1, \Om_2$, $\cdots$, $\Om_{N_s}$, then partition the time interval $[0,T]$ into $N_t$ subintervals using the partition: $0<T_1<T_2<\cdots<T_{N_t}$.
By coupling all the space subdomains and time subintervals, a decomposition of $\Theta$ is $\Theta = \cup_{i=1}^{N_s}
(\cup_{j=1}^{N_t} \Theta_{ij})$, where $\Theta_{ij} = \Om_i \times (T_{j-1}, T_{j})$.
For convenience, the number of subdomains, i.e. $N_sN_t$, is equal to the number of processors.
These subdomains $\Theta_{ij}$ are then extended to $\Theta'_{ij}$ to overlap each other.
The boundary of each subdomain is extended by an integral number of mesh cells in each dimension,
and we trim the cells outside of $\Theta$.
The corresponding overlapping decomposition of $\Theta$ is $\Theta = \cup_{i=1}^{N_s}
(\cup_{j=1}^{N_t} \Theta'_{ij})$.
See the left figure of Figure \ref{fig:ddm} for the overlapping extension.

The matrix on each subdomain $\Theta_{ij}'=\Om_{i}'\times (T'_{j-1}, T'_{j})$, $i=1,2$, $\cdots$, $N_s$, $j=1,2$, $\cdots$, $N_t$ is the discretized version of the following system of PDEs
\begin{equation}
\begin{cases}
\displaystyle\frac{\partial C}{\partial t}=\nabla\cdot(a(\mathbf{x})\nabla C)-\nabla\cdot (\mathbf{v}(\mathbf{x})C)+f(\x, t)\,, ~ (\x,t)\in \Theta_{ij}' \\
\displaystyle\frac{\partial G}{\partial t}=-\nabla \cdot (a(\mathbf{x})\nabla G)-\mathbf{v}(\mathbf{x})\cdot \nabla G \\
{\,\ } \hskip1truecm +A(\x)(C(\x,t)-C^{\varepsilon}(\x, t))\,, ~ ~ (\x,t)\in \Theta_{ij}' \\
\beta_1\displaystyle\frac{\partial^2 f}{\partial t^2}+\beta_2\Delta f+G = 0\,, ~ (\x,t)\in \Theta_{ij}'
\end{cases}
\label{eq:subprob}
\end{equation}
with the following boundary conditions
\begin{eqnarray}\label{eq:subprobBCs}
C(\mathbf{x},t)=0, \q G(\mathbf{x},t)=0,\q f(\mathbf{x},t)=0\,,~ \x\in \partial\Om_{i}'\,,~t\in [T'_{j-1}, T'_{j}]
\end{eqnarray}
along with the initial and terminal time boundary conditions
\begin{equation}\label{eq:subprobIC}
\begin{cases}
C(\mathbf{x},T'_{j-1})=0, \q G(\x, T'_{j-1})=0, \q f(\x, T'_{j-1})=0\,,~ \x\in \partial\Om_{i}'\\
C(\mathbf{x},T'_{j})=0, \q G(\x, T'_{j})=0, \q f(\x, T'_{j})=0\,,~ \x\in \partial\Om_{i}'\,.
\end{cases}
\end{equation}

One may notice from (\ref{eq:subprobIC}) that the homogenous Dirichlet boundary conditions are applied
in each time interval $(T'_{j-1}, T'_{j})$, so the solution of each subdomain problem is not really physical.
This is one of the major differences between the space-time Schwarz method and
the parareal algorithm \cite{LMT01}.
The time boundary condition for each subproblem of the parareal algorithm
is obtained by an interpolation of the coarse solution, and if the coarse mesh is fine enough, the solution of the subdomain problem is physical. As a result, the parareal algorithms
can be used as a solver, but our space-time Schwarz method can only be used as a preconditioner.
Surprisingly, as we shall see from our numerical experiments in Section~\ref{sec:Numerical Examples},
the Schwarz method is an excellent preconditioner even though the time boundary conditions violate the physics.

We solve the subdomain problems the same as for the global problem (\ref{eq:kktsys2}),
no time-marching is performed in our new algorithm, and
all unknowns affiliated with each subdomain are solved simultaneously.
 Let $M_{ij}$ be the matrix generated in the same way as the global matrix $F_h$ in (\ref{eq:kktsys2})
 but for the subproblem (\ref{eq:subprob})-(\ref{eq:subprobIC}),
 and $\tilde{M}_{ij}^{-1}$ be an exact or approximate inverse of $M_{ij}$.
 Denoting the restriction matrix from $\Theta$ to the subdomain $\Theta'_{ij}$ by $R^{\delta}_{ij}$, with overlapping size $\delta$, the space-time restricted Schwarz preconditioner \cite{CS99}  can be now formulated
 as
\[{\,\ } \hskip1truecm M_{one-level}^{-1} = \sum_{j=1}^{N_t}\sum_{i=1}^{N_s} (R^{\delta}_{ij})^T\tilde{M}_{ij}^{-1}R^{0}_{ij}\,.\]

As it is well known, any one-level domain decomposition methods are not scalable with the increasing number
of subdomains or processors.  Instead one should have multilevel methods in order to observe
possible scalable effects \cite{ACM99,SBG04}.
We now propose a two-level space-time additive Schwarz preconditioner.
To do so, we partition $\Om$ with a fine mesh $\Om^h$ and a coarse mesh $\Om^c$.
For the time interval, we have a fine partition $P^{\tau}$ and a coarse partition $P^{\tau_c}$  with $\tau < \tau_c$.
We will adopt a nested mesh, i.e., the nodal points of the coarse mesh $\Om^c\times P^{\tau_c}$
are a subset of the nodal points of the fine mesh $\Om^h\times P^{\tau}$.
In practice, the size of the coarse mesh should be adjusted properly to obtain the best performance.
On the fine level, we simply apply the previously defined one-level space-time additive Schwarz preconditioner; and to efficiently solve the coarse problem, a parallel coarse preconditioner is also necessary.
Here we use the overlapping space-time additive Schwarz preconditioner and for simplicity divide $\Om^c\times P^{\tau_c}$ into the same number of subdomains as on the fine level, using the non-overlapping decomposition $\Theta = \cup_{i=1}^{N_s}
(\cup_{j=1}^{N_t} \Theta_{ij})$.
When the subdomains are extended to overlapping ones, the overlapping size is not necessarily the same as that on the fine mesh.
See the right figure of Figure \ref{fig:ddm} for a coarse version of the space-time decomposition.
We denote the preconditioner for the coarse level by $M^{-1}_c$, which is defined by
\[{\,\ } \hskip1truecm M_{c}^{-1} = \sum_{j=1}^{N_t}\sum_{i=1}^{N_s} (R^{\delta_c}_{ij,c})^T\tilde{M}_{ij,c}^{-1}R^{0}_{ij,c}\,,\]
where $\delta_c$ is the overlapping size on the coarse mesh.
Here the matrix $\tilde{M}^{-1}_{ij,c}$ is an approximate inverse of $M_{ij,c}$ which is obtained by a discretization of (\ref{eq:subprob})-(\ref{eq:subprobIC}) on the coarse mesh on $\Theta'_{ij}$.

To combine the coarse preconditioner with the fine mesh preconditioner, we need a restriction operator $I^c_h$
from the fine to coarse mesh and an interpolation operator $I^h_c$ from the coarse to fine mesh.
For our currently used nested structured mesh and linear finite elements,
$I^h_c$ is easily obtained using a linear interpolation on the coarse mesh and $I^c_h = (I^h_c)^T$.
We note that when the coarse and fine meshes are nested, instead of using $I^c_h = (I^h_c)^T$,
we may take $I^c_h$ to be a simple restriction, e.g., the identity one which assigns the values on the coarse mesh
using the same values on the fine mesh.
In general, the coarse preconditioner and the fine preconditioner can be combined additively or multiplicatively.
According to our experiments, the following multiplicative version works well:
\begin{equation}
\begin{cases}
y=I^h_c F_{c}^{-1} I^c_h x\\
M^{-1}_{two-level}x = y+M_{one-level}^{-1}(x-F_h y)\,,
\end{cases}
\end{equation}
where $F_c^{-1}$ corresponds to the GMRES solver right-preconditioned by $M_{c}^{-1}$ on the coarse level,
and $F_h$ is the discrete KKT system (\ref{eq:kktsys2})  on the fine level.

\section{Numerical experiments}\label{sec:Numerical Examples}
In this section we present some numerical experiments to study the parallel performance and robustness of the newly developed algorithms.
When using the one-level preconditioner, we use a restarted GMRES method (restarts at 50) to solve the preconditioned system;
when using the two-level preconditioner,  we use the restarted flexible GMRES (fGMRES) method \cite{S93} (restarts at 30),
considering the fact that
the overall preconditioner changes from iteration to iteration because of the iterative coarse solver.
Although fGMRES needs more memory than GMRES, we have observed
its number of iterations can be significantly reduced.
The relative convergence tolerance of both GMRES and fGMRES is set to be $10^{-6}$.
The initial guesses for both GMRES and fGMRES method are zero.
The size of the overlap between two neighbouring subdomains, denoted by $iovlp$, is set to be 1 unless otherwise specified.
The subsystem on each subdomain is solved by an incomplete LU factorization ILU($k$), with $k$ being
its fill-in level, and $k=0$ if not specified.
The algorithms are implemented based on the Portable, Extensible Toolkit for Scientific computation (PETSc) \cite{BBEG10}
using run on a Dawning TC3600 blade server system at the National Supercomputing Center in Shenzhen, China with a
1.271 PFlops/s peak performance.

In our computations, the settings for the model system (\ref{eq:cfeq}) are taken as follows.
The computational domain, the terminal time and the initial condition are taken to be $\Om=(-2,2)^3$,
$T=1$ and $C(\cdot,0)=0$ respectively. Let $L=S=H=2$, then
the homogeneous Dirichlet and Neumann conditions in (\ref{eq:cfeq}) are respectively imposed
on $\Gamma_1=\{\mathbf{x}=(x_1,x_2,x_3); ~|x_1|=L~\mbox{or }~|x_2|=S\}$  and  $\Gamma_2=\{\mathbf{x}=(x_1,x_2,x_3);  ~|x_3|=H\}$.
Furthermore, the diffusivity and convective coefficients are set to be
$a(\mathbf{x})=1.0$ and $\mathbf{v}(\mathbf{x})=(1.0,1.0,1.0)^T$.

In order to generate the observation data, we solve the forward convection-diffusion equation
(\ref{eq:cfeq}) on a very fine mesh with a small time step size, and the resulting approximate
solution $C(\x,t)$ is used as the noise-free observation data.
Then a random noise is added in the following form at the locations where the measurements are taken:
\begin{equation*}
C^{\epsilon}(\mathbf{x}_i,t)=C(\mathbf{x}_i,t)+\epsilon\, r\,C(\mathbf{x}_i,t), \quad i=1,\cdots,s\,.
\end{equation*}
Here $r$ is a random function with the standard Gaussian distribution, and $\epsilon$ is the noise level.
In our numerical experiments, $\epsilon = 1\%$  if not specified otherwise.

The numerical tests are designed to investigate the reconstruction effects with different types
of three-dimensional sources
by the proposed one- or two-level space-time Schwarz method, as well as
the robustness of the algorithm with respect to different noise level, different regularizations and
different amount of measurement data. In addition, parallel efficiency of the proposed algorithms is
also studied.

\subsection{Reconstruction of 3D sources}
We devote this subsection to test the numerical reconstruction of three representative 3D sources
by the proposed one-level space-time method, with $np=256$ processors. Each of the three examples are constructed
with its own special difficulty.

\begin{figure}
\centering{
\includegraphics[width=6cm]{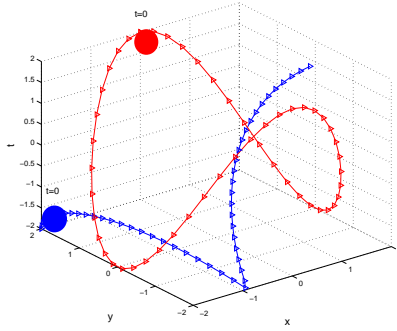}
\vspace{-0.2in}
\caption{ The traces of two moving sources.}
\label{fig:ex2trace}
}
\end{figure}

\textbf{Example 1: two Gaussian sources}.
This example tests
two moving Gaussian sources in $\Om$, namely the source $f$ takes the form:
\[f(\x,t)=\displaystyle\sum_{i=1}^2 \mbox{exp}\Big(-\displaystyle\frac{(x-x_i)^2+(y-y_i)^2+(z-z_i)^2}{a^2} \Big)\,,\]
with $a=2.0$ and two moving centers of the sources are given by
\begin{equation}
\begin{cases}
(x_1,y_1,z_1)=(L\sin(2\pi t),S\cos(2\pi t),H\cos(4\pi t))\, \\
(x_2,y_2,z_2)=(L-2L|\cos(4t)|, -S+2S|\cos(4t)|, -H+2Ht^2) \,.
\end{cases}
\end{equation}
The moving traces of the sources are shown in Figure \ref{fig:ex2trace}.

\begin{figure}
\centering
\includegraphics[width=1.0\textwidth]{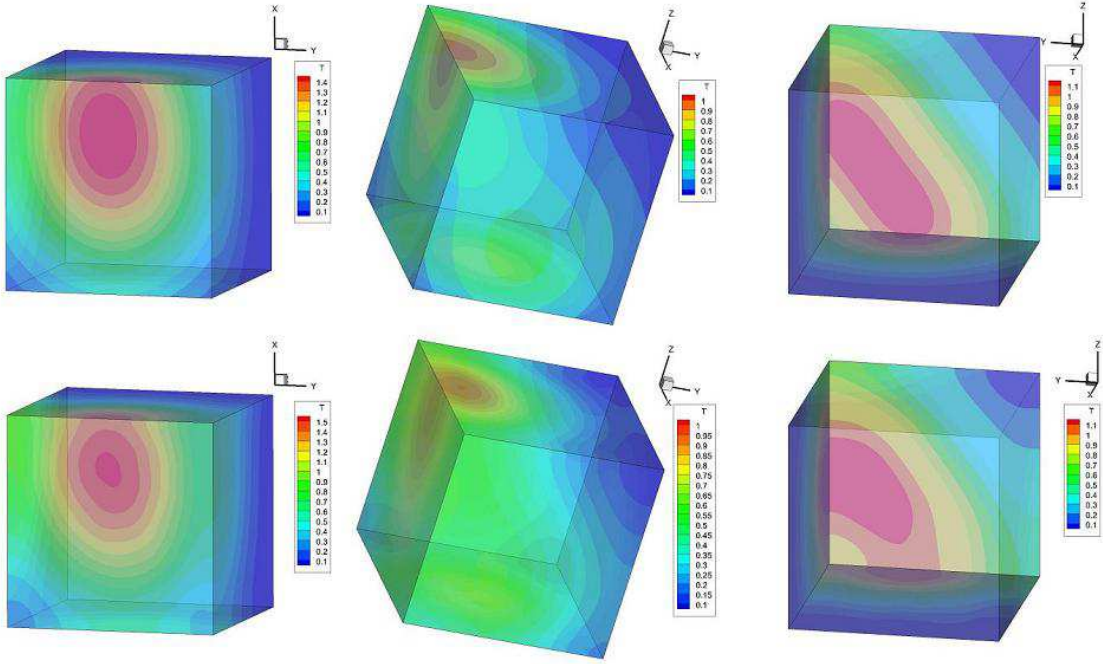}
\caption{ Example 1: the source reconstructions at three moments $t=10/39, 20/39, 30/39$ with
measurements collected at the mesh $14\times14\times14$ (bottom),
comparable with the exact source distribution (top).}\label{fig:ex1}
\end{figure}

\begin{figure}
\centering
\includegraphics[width=1.0\textwidth]{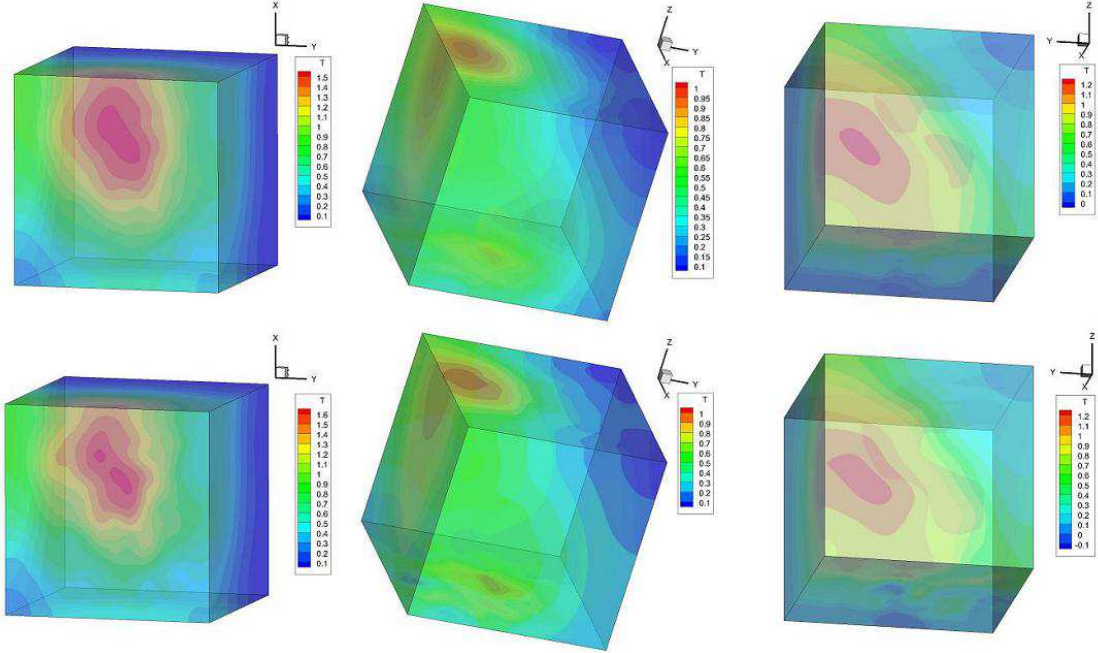}
\vspace{-0.3in}
\caption{ Example 1: source reconstructions with noise level $\epsilon = 5\%$ (top) and $\epsilon = 10\%$ (bottom).}\label{fig:ex1noise}
\end{figure}

In the first experiment, we use the mesh $40\times40\times40$ and
the time step size of $1/39$ for the inversion process.
And the measurements are taken on the mesh
$14\times 14\times 14$, which is uniformly located in $\Om$.
The regularization parameters are set to be
$\beta_1=3.6\times 10^{-5}$ and $\beta_2=3.6\times 10^{-3}$.
In Figure \ref{fig:ex1}, the numerically reconstructed sources are compared
with the exact one at three moments $t=10/39, 20/39, 30/39$. We can see
that the source locations and intensities are quite close to the true values at three chosen moments.
Then we increase the noise level to $\epsilon = 5\%$ and $\epsilon = 10\%$, still with
the same set of parameters. The reconstruction results are shown in Figure \ref{fig:ex1noise}.
We can observe that the reconstructed profiles deteriorate and become oscillatory as the noise level increases.
This is naturally expected since the ill-posedness of the inverse source problem increases
with the noise level.

\textbf{Example 2: Four constant sources}.
Appropriate choices of regularizations are important for the inversion process.
In the previous example we have used a $H^1$-$H^1$ Tikhonov regularization in both space and time.
In this example, we intend to compare the $H^1$-$H^1$ regularization
with the following $H^1$-$L^2$ regularization
\begin{equation*}
 \tilde{N}_{\beta}(f)=\displaystyle\frac{\beta_1}{2}\int_0^T\int_{\Om} | \dot{f}(\x, t)|^2 d\x dt + \displaystyle\frac{\beta_2}{2}\int_0^T\int_{\Om} f^2 d\x  dt\,.
\end{equation*}
For the comparisons, we consider the case
in which four constant sources move along the diagonals of the cube to their far corner.
The four sources are specified by
\begin{equation*}
f_i(\x,t)=a_i~~\mbox{for }~|x-x_i|<0.4, ~~|y-y_i|<0.4, ~~|z-z_i|<0.4
\end{equation*}
where $a_1 = a_4 = 2.0$, $a_2 = a_3 = 1.0$,  and the traces of the four sources are
described by
\begin{equation*}
\begin{cases}
(x_1,y_1,z_1) = (-L+2Lt,-S+2St,H-2Ht)\\
(x_2,y_2,z_2) = (L-2Lt,S-2St,-H+2Ht)\\
(x_3,y_3,z_3) = (L-2Lt,-S+2St,-H+2Ht)\\
(x_4,y_4,z_4) = (-L+2Lt, S-2St, H-2Ht)\,.
\end{cases}
\end{equation*}
Same mesh and measurements are used as in Example 1, and
the regularization parameters are set to be
$\beta_1= 10^{-5}, \beta_2= 10^{-3}$ in $N_{\beta}(f)$, and $\beta_1= 10^{-5}, \beta_2= 10^{-8}$
in $\tilde{N}_{\beta}(f)$, respectively.
The reconstruction results are compared with the true solution
at three moments $t=10/39, 20/39, 30/39$, and
two slices at $x=0.95$ and $x=-0.95$.
\begin{figure}
\centering
\includegraphics[width=1.0\textwidth]{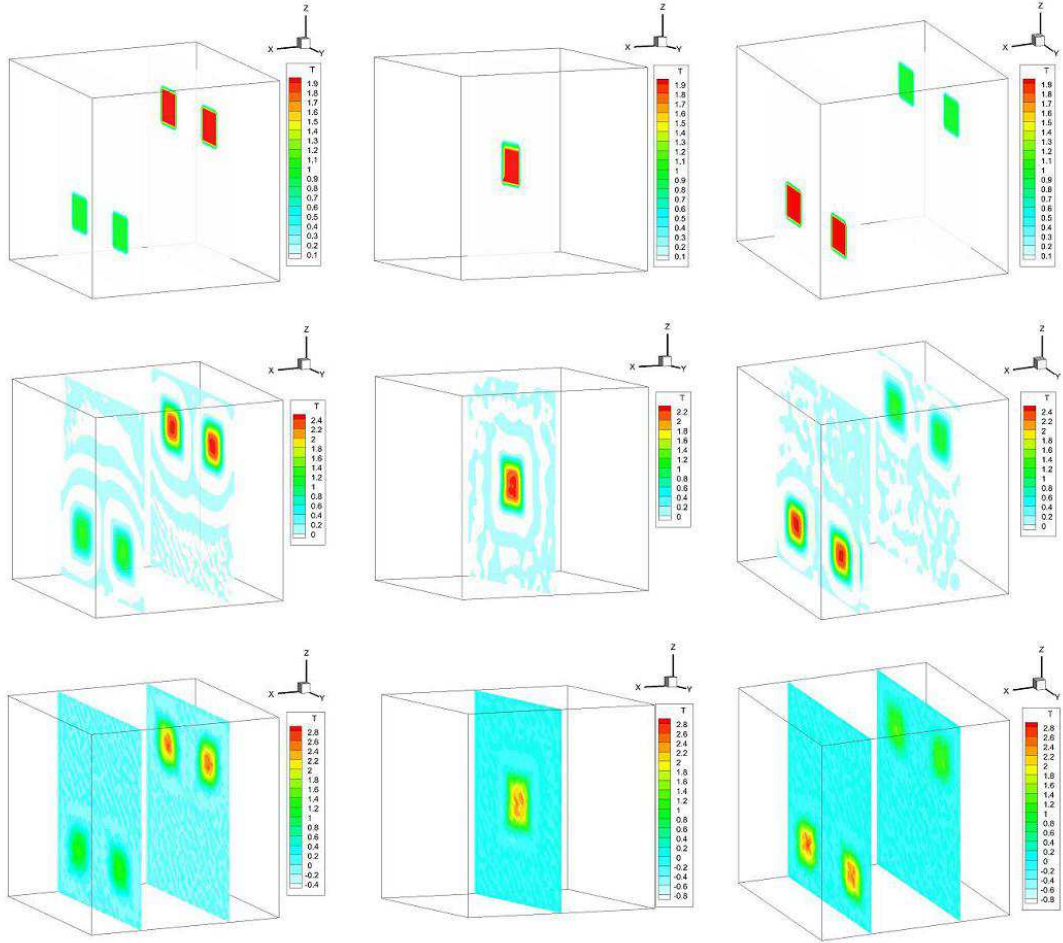}
\vspace{-0.3in}
\caption{ Example 2: the source reconstructions with $H^1$-$H^1$ regularization (mid) and $H^1$-$L^2$ regularization (bottom), compared with the exact solution (top).}\label{fig:ex2}
\end{figure}
It is observed from Figure \ref{fig:ex2} that the resolution of the source profile is much better with the
$H^1$-$H^1$ regularization $N_{\beta}(f)$ than with the $H^1$-$L^2$ regularization $\tilde{N}_{\beta}(f)$,
and the latter presents a reconstruction process that is much less stable and much more oscillatory.

\begin{figure}
\centering
\includegraphics[width=1.0\textwidth]{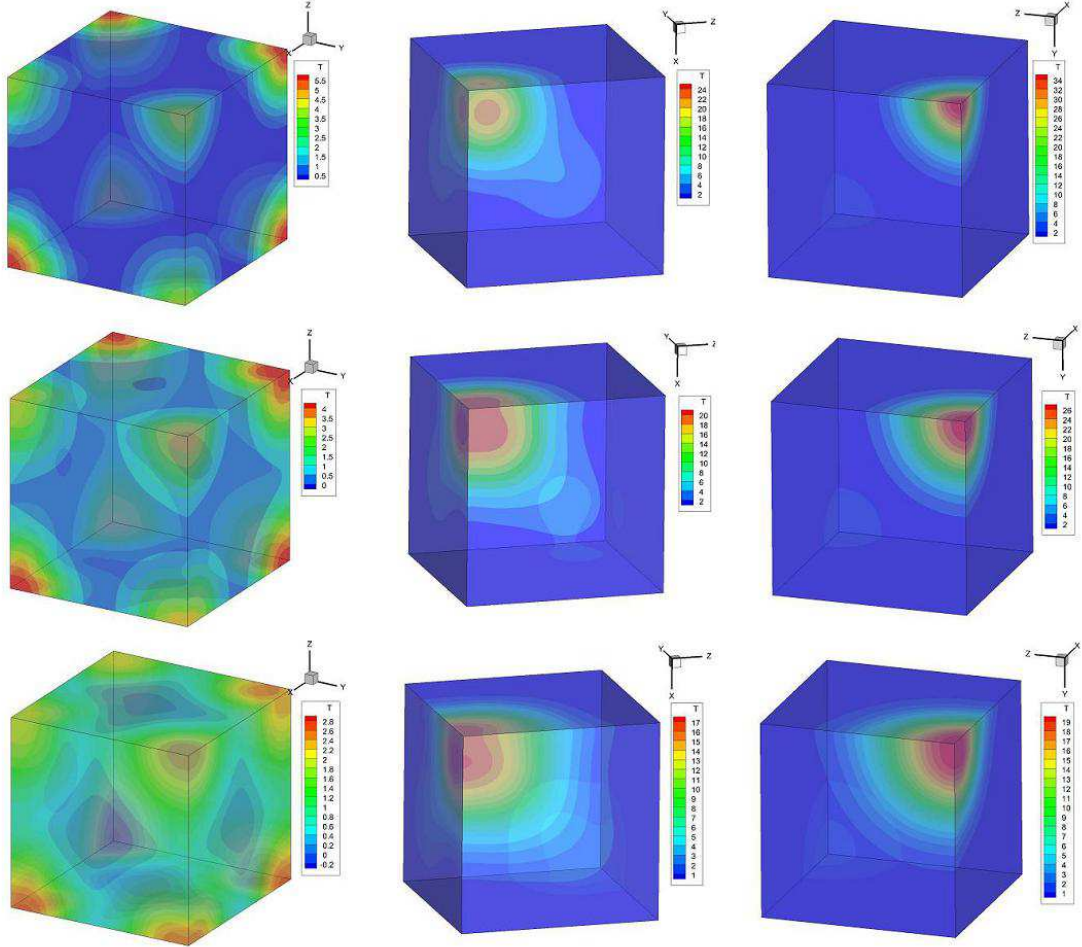}
\vspace{-0.3in}
\caption{ Example 3: the source reconstructions with measurements collected at the mesh $21\times21\times21$
(mid) and $9\times9\times9$ (bottom), compared with the exact solution (top)}\label{fig:ex3}
\end{figure}

\textbf{Example 3: Eight moving  sources}.
This last example presents a very challenging case
that eight Gaussian sources are initially located at the corners of the physical cubic domain,
then move inside the cube following their own traces given below. The Gaussian sources are described by
\begin{equation*}
f(\x,t)=\displaystyle\sum_{i=1}^8 a_i e^{-(x-x_i)^2+(y-y_i)^2+(z-z_i)^2}\,,
\end{equation*}
where the coefficients $a_i$ and the source traces are represented by
\begin{align*}
a_1=a_2=&a_3=a_4=4.0; ~a_5=a_6=a_7=a_8=6.0,\\
(x_1,y_1,z_1)&=\left(-L+2L(1-t), -S+2S(1-t), -H+2H(1-t)\right)\, \\
(x_2,y_2,z_2)&=\left(-L+2Lt, -S+2St, -H+2Ht\right) \,\\
(x_3,y_3,z_3)&=\left(-L+2L\cos(\pi t)^2(1-t), -S+2S\sin(\pi t)^2 t, -H+2H\cos(\pi t)^2(1-t)\right)\, \\
(x_4,y_4,z_4)&=\left(-L+2L\cos(\pi t)^2(1-t), -S+2S\cos(\pi t)^2(1-t), -H+2H\sin(\pi t)^2 t)\right) \,\\
(x_5,y_5,z_5)&=\left(-L+2L\cos(2\pi t)^2\cos\left(\pi/2 t\right), -S+2S\sin(\pi \right. t)^2\sin\left(\pi/2 t\right), \\
&\left. -H+2H\sin(\pi t)^2\sin\left(\pi/2 t\right)\right)\, \\
(x_6,y_6,z_6)&=\left(-L+2L\sin(\pi t)^2\sin\left(\pi/2 t\right), -S+2S\cos(2\pi \right. t)^2\cos\left(\pi/2 t\right), \\
&\left. -H+2H\sin(\pi t)^2\sin\left(\pi/2 t\right)\right)\\
(x_7,y_7,z_7)&=\left(-L+2L\sin(\pi t)^2\sin\left(\pi/2 t\right), -S+2S\sin(\pi \right. t)^2\sin\left(\pi/2 t\right), \\
&\left. -H+2H\cos(2\pi t)^2\cos\left(\pi/2 t\right)\right)\,\\
(x_8,y_8,z_8)&=\left(-L+2L\sin(\pi t)^2\sin\left(\pi/2 t\right), -S+2S\cos(2\pi\right. t)^2\cos\left(\pi/2 t\right), \\
&\left. -H+2H\cos(2\pi t)^2\cos\left(\pi/2 t\right)\right) \,.
\end{align*}
We shall use the mesh $64\times64\times64$ and the time step size $1/47$, with
two regularization parameters $\beta_1=3.6\times 10^{-5}$ and $\beta_2=3.6\times 10^{-1}$.
We compare the results recovered by two sets of measurements, collected at
two meshes $21\times21\times21$ and
$9\times9\times9$ respectively, which are both uniformly distributed in $\Om$,
with the exact solution shown in Figure \ref{fig:ex3} (top),
at three time moments $t= 0.0, 10/47, 1.0$.
Clearly better reconstructions are observed for the case with more measurements collected at the finer mesh
$21\times21\times21$,
though the coarser mesh $9\times9\times9$ is good enough for locating the sources, only with
their recovered source intensities smaller than the true values.

\subsection{Performance in parallel efficiency}
In the previous subsection, we have shown with 3 representative examples that the proposed algorithm
can successfully recover the intensities and distributions of unsteady sources and is
robust with respect to the noise in the data, the choice of Tikhonov regularizations and the number of measurements.
These numerical simulations are all computed using the proposed
one-level space-time method with $np=256$ processors.
In this section, we focus on our proposed two-level space-time method and study its parallel efficiency
with respect to the number of ILU fill-in levels, namely the number $k$ in ILU($k$), and the overlap size $iovlp$.
We also compare the number of iterations and the total compute time of the one-level and two-level methods with increasing degrees of freedoms (DOFs) and the number of processors.

First we will test how the number of fGMRES iterations and the total compute time of the two-level method
change with different ILU fill-in levels.
We use the coarse mesh $21\times21\times21$ with the time step $1/20$,
and the fine mesh $41\times41\times41$ with the time step $1/40$
for Example 1, 2, and 3, and the overlap size $iovlp=1$. We see that
the total number of degrees of freedom on the fine mesh is 16 times of the one
on the coarse mesh.
Table \ref{tab:ilulevel} shows the comparison with $np= 256$  processors. Column 2-3, 4-5 and 6-7
present the results for Example 1,  2 and 3 respectively.
It is observed that as the fill-in level increases the number of fGMRES iterations decreases, but the total compute time increases. When the fill-in level increases to 3, the compute time increases significantly and the number of iterations only reduces
by 3 times. This suggests a suitable fill-in level to be $ilulevel =0$ or 1.
\begin{table}
\centering
\caption{Effects of ILU fill-in levels on the two-level method for Example 1 (columns 2-3),
Example 2 (columns 4-5), and Example 3 (columns 6-7).}
\label{tab:ilulevel}
\begin{tabular}{lllllll}
\hline\noalign{\smallskip}
ILU($k$)&Its&Time (sec)&Its&Time (sec)&Its&Time (sec)\\
\noalign{\smallskip}\hline\noalign{\smallskip}
0&47&10.498&55&12.448&81&17.238\\
1&28&33.633&36&47.766&60&49.622\\
2&18&230.552&23&232.914&48&257.798\\
3&15&1121.469&20&1132.841&45&1165.203\\
\noalign{\smallskip}\hline
\end{tabular}
\end{table}

Next we look at the impact of the overlap size.
We still use the same fine and coarse meshes for all examples, and ILU(0) for
the solver for each subdomain problem on both the coarse and fine meshes.
The overlap size on the coarse mesh is set to be 1.
We test different overlap sizes on the fine level, and the results are given in Table \ref{tab:overlap}.
It is observed that when the overlap size increases from 1 to 2 and then to 4, the number of fGMRES iterations decreases slowly and the total compute time increases.
So we shall mostly use $iovlp =1$ in our subsequent computations.
\begin{table}
\centering
\caption{Effects of the overlap size on the two-level method for
Example 1 (columns 2-3),  Example 2 (columns 4-5), and Example 3 (columns 6-7).}
\label{tab:overlap}
\begin{tabular}{lllllll}
\hline\noalign{\smallskip}
$iovlp$&Its&Time (sec)&Its&Time (sec)&Its&Time (sec)\\
\noalign{\smallskip}\hline\noalign{\smallskip}
1&47&10.498&55&12.448&81&17.238\\
2&39&13.071&51&23.663&69&27.952\\
4&37&27.423&49&45.225&68&47.032\\
\noalign{\smallskip}\hline
\end{tabular}
\end{table}

Lastly we compare the performance of the one-level and two-level space-time Schwarz preconditioners in Tables \ref{tab:twolevelmesh} and \ref{tab:twolevelnp}.
On the coarse level, a restarted GMRES is used, with the one-level space-time Schwarz preconditioner.
ILU(0) is used as the local preconditioner on each subdomain and the coarse overlap size is set to be 1.
A tighter convergence tolerance on the coarse mesh can reduce the number of outer fGMRES iterations, but often increases the total compute time.
In the following numerical examples, we set the tolerance to be $10^{-1}$ and the maximum number of GMRES iterations
to 4 on the coarse mesh.
Moreover, the mesh size of the coarse mesh is also an important factor for the performance.
If the mesh is too coarse, both the number of outer iterations and the total compute time increase;
on the other hand, if the mesh is not coarse enough, too much time is spent for the coarse solver, the number of outer iterations may decrease significantly, but the compute time may increase.

In the following experiments for Example 1, 2 and 3,  we use three sets of  fine meshes, $33\times33\times33$, $49\times49\times49$ and $67\times67\times67$, and the corresponding time steps are $1/32$, $1/48$ and $1/66$ respectively, while
the coarse meshes are chosen to be $17\times17\times17$, $17\times17\times17$ and $23\times23\times23$,
with the corresponding time steps being $1/16$, $1/48$ and $1/66$.
So the DOFs on the fine meshes are 16, 27 and 27 times of the ones on the coarse meshes
for Example 1, 2 and 3 respectively. We use $np=64, 128$ and 512 processors for the three sets of meshes respectively and compare their performance with the one-level method in Table  \ref{tab:twolevelmesh}.
Savings in terms of the number of iterations and the total compute time are obtained for the two-level method with all three sets of meshes. As we observe that the number of iterations of the two-level method is mostly reduced by at least
4 times compared to the one for the one-level method, but the compute time is usually reduced by 2 to 4 times.

Next we fix the space mesh to be $49\times49\times49$ and the time step to be $1/48$,
resulting in a very large-scale discrete system with 17,294,403 DOFs.
For the two-level method, we set the coarse mesh to be $17\times17\times17$ with the time step $1/48$,
which implies that the DOFs on the fine mesh is about 27 times of the ones on the coarse mesh.
Then the problem is solved with $np = 128, 256, 512$, and 1024 processors respectively.
The performance results of the one-level and two-level methods are presented in Table  \ref{tab:twolevelnp}.
We observe that when the number of sources is small, both the one-level and two-level methods are scalable with up to 512 processors, but the two-level method takes much less compute time.
The strong scalability deteriorates when the number of processors is too large for the size of the problems.
 As the number of sources increases, the scalability becomes slightly worse
 for both one-level and two-level methods, even though the two-level method is still faster in terms of the total compute time.
\begin{table}
\centering
\caption{Comparisons between the one-level and two-level space-time preconditioners
for Examples 1-3 with different meshes.}
\label{tab:twolevelmesh}
\begin{tabular}{llllll}
\hline\noalign{\smallskip}
\multicolumn{6}{c}{Ex1}\\
\noalign{\smallskip}\hline\noalign{\smallskip}
$np$&Mesh&$M$&$level$&Its&Time (sec)\\
\noalign{\smallskip}\hline
64&$33\times33\times33$&33&1&175&53.635\\
&&&2&57&20.653\\
\noalign{\smallskip}\hline
128&$49\times49\times49$&49&1&346&200.664\\
&&&2&83&47.812\\
\noalign{\smallskip}\hline
512&$67\times67\times67$&67&1&491&675.985\\
&&&2&105&212.72\\
\noalign{\smallskip}\hline
\multicolumn{6}{c}{Ex2}\\
\noalign{\smallskip}\hline
$np$&Mesh&$M$&$level$&Its&Time (sec)\\
\noalign{\smallskip}\hline
64&$33\times33\times33$&33&1&228&72.338\\
&&&2&77&20.246\\
\noalign{\smallskip}\hline
128&$49\times49\times49$&49&1&365&214.058\\
&&&2&85&47.078\\
\noalign{\smallskip}\hline
512&$67\times67\times67$&67&1&599&841.652\\
&&&2&121&216.92\\
\noalign{\smallskip}\hline
\multicolumn{6}{c}{Ex3}\\
\noalign{\smallskip}\hline
$np$&Mesh&$M$&$level$&Its&Time (sec)\\
\noalign{\smallskip}\hline
64&$33\times33\times33$&33&1&297&82.834\\
&&&2&76&21.738\\
\noalign{\smallskip}\hline
128&$49\times49\times49$&49&1&405&238.712\\
&&&2&93&57.244\\
\noalign{\smallskip}\hline
512&$67\times67\times67$&67&1&716&872.766\\
&&&2&137&263.222\\
\noalign{\smallskip}\hline
\end{tabular}
\end{table}

\begin{table}
\centering
\caption{Comparisons between the one-level and two-level space-time preconditioners
for Examples 1-3 with different number of processors.}
\label{tab:twolevelnp}
\begin{tabular}{llllllll}
\hline\noalign{\smallskip}
&&Ex1&&Ex2&&Ex3&\\
\noalign{\smallskip}\hline\noalign{\smallskip}
$np$&$level$&Its&Time (sec)&Its&Time (sec)&Its&Time (sec)\\
\noalign{\smallskip}\hline
128&1&346&200.664&365&214.815&405&238.712\\
&2&83&47.812&85&47.072&93&57.244\\
\noalign{\smallskip}\hline
256&1&343&127.035&363&152.334&408&145.213\\
&2&82&24.744&87&26.424&90&36.307\\
\noalign{\smallskip}\hline
512&1&343&69.482&363&95.707&400&101.343\\
&2&82&16.461&101&19.453&100&18.611\\
\noalign{\smallskip}\hline
1024&1&351&41.821&393&58.785&433&59.534\\
&2&85&10.132&100&11.352&104&15.815\\
\noalign{\smallskip}\hline
\end{tabular}
\end{table}

\section{Concluding remarks}\label{sec:conclusions}
In this work we have proposed and studied a new fully implicit, space-time coupled, mixed finite element and finite difference discretization method, and a parallel one- and two-level domain decomposition solver for
the three-dimensional unsteady inverse convection-diffusion problem.
With a suitable number of measurements,
this all-at-once approach provides acceptable reconstruction of the physical sources in space and time simultaneously.
The classical overlapping Schwarz preconditioner is extended successfully to the coupled space-time problem with a homogenous Dirichlet boundary condition applied on both the spatial and temporal part of the space-time
subdomain boundaries.
The one-level method is easier to implement, but the two-level hybrid space-time Schwarz method performs much better in terms of the number of iterations and the total compute time.
Good scalability results were obtained for problems with more than 17 millions degrees of freedom on a supercomputer with more than 1,000 processors.
The approach is promising to more general unsteady inverse problems in large-scale applications.

\appendix
\section{The discrete structure of the KKT system}
The KKT system (\ref{eq:kkt3})-(\ref{eq:kkt6}) is formulated as follows:
\begin{equation}
\begin{cases}\label{eq:kktsys}
(\partial_{\tau} C_h^n, v_h) + (a\nabla  \bar{C}_h^n, \nabla v_h) + (\nabla\cdot (\mathbf{v}\bar{C}_h^n),v_h) =  (\bar{f}^n_h, v_h)+\langle\bar{q}^n,v_h\rangle_{\Gamma_2}, ~~\forall\,v_h\in \mathring{V}^h\,\\
-(\partial_{\tau} G_h^n, w_h) + (a\nabla  \bar{G}_h^n, \nabla w_h) + (\nabla\cdot (\mathbf{v}w_h),\bar{G}_h^n)
\\
=-(A(\x)(\bar{C}_h^{n}(\x,t)-\bar{C}^{\epsilon,n}(\x, t)), w_h), ~~\forall w_h\in \mathring{V}^h \,\\
-  (G_h^n, g_h^{\tau}) +\beta_1( \partial_{\tau} f_h^n,\partial_{\tau} g_h^{\tau}) +  \beta_2 (\nabla f_h^n, \nabla  g_h^{\tau})=0,  ~~\forall\,g_h^{\tau}\in W_h^{\tau}\,.
\end{cases}
\end{equation}
To better understand the discrete structure of (\ref{eq:kktsys}), we denote the identity and zero matrices as $I$ and $\mathbf{0}$ respectively, and the basis functions of the finite element spaces $V^h$ and $W_h^{\tau}$ by $\phi=(\phi_i)^T$, $i=1$, $\cdots$, $N$ and $g_j^n$, $j=1$, $\cdots$, $N$, $n=0$, $\cdots$, $M$, respectively, let
\begin{align*}
&A=(a_{ij})_{i,j=1,\cdots,N},\quad a_{ij}=(a \nabla \phi_i, \nabla \phi_j)\\
&B=(b_{ij})_{i,j=1,\cdots,N},\quad b_{ij}=( \phi_i,  \phi_j)\\
&E=(e_{ij})_{i,j=1,\cdots,N},\quad e_{ij}=(\nabla\cdot(\mathbf{v} \phi_i), \phi_j)\\
&L^{mn}=(l_{ij}^{mn})_{i,j=1,\cdots,N, 0\leq m, n\leq M},\quad l_{ij}^{mn}=\left(\displaystyle\frac{\partial g_i^m}{\partial t}, \displaystyle\frac{\partial g_j^n}{\partial t}\right)\\
&K^{mn}=(k_{ij}^{mn})_{i,j=1,\cdots,N, 0\leq m, n\leq M},\quad k_{ij}^{mn}=( \nabla g_i^m, \nabla g_j^n)\\
&D^{mn}= (d_{ij}^{mn})_{i,j=1,\cdots,N, 0\leq m, n\leq M},\quad d_{ij}^{mn}=( g_i^m,  g_j^n)\,,
\end{align*}
and based on these element matrices we define
\begin{align*}
&A_1=B+\displaystyle\frac{\tau}{2}(A+E),\quad A_2=-B+\displaystyle\frac{\tau}{2}(A+E)\quad \quad \quad \\
&B_1=B+\displaystyle\frac{\tau}{2}(A+E^T),\quad B_2 =-B+\displaystyle\frac{\tau}{2}(A+E^T)\quad \quad \quad \\
&B_3= \mbox{zeros except 1 at the measurement locations}\quad \quad \quad \\
&W^{mn}=\beta_1 L^{mn}+ \beta_2 K^{mn}\,,\quad \quad \quad
\end{align*}
Then the system (\ref{eq:kktsys}) takes the following form
\begin{equation*}
\left(\begin{array}{ccc}BC&BG&Bf\end{array}\right)\left(\begin{array}{c}C^0\\C^1\\ \vdots\\C^{M-2}\\C^{M-1}\\C^{M}\\G^0\\G^1\\G^2\\ \vdots \\G^{M-2}\\ G^{M-1}\\G^{M}\\f^0\\f^1\\f^2\\ \vdots \\f^{M-2}\\f^{M-1}\\f^{M}\\ \end{array}\right)
=\left(\begin{array}{c}C^0\\ \langle\bar{q}^1,\phi\rangle_{\Gamma_2} \\ \vdots\\ \langle\bar{q}^{M-1},\phi\rangle_{\Gamma_2} \\ \langle\bar{q}^{M},\phi\rangle_{\Gamma_2}\\ \tau/2 B_3 (C^{\epsilon,0}+ C^{\epsilon,1})\\ \vdots \\  \tau/2 B_3 (C^{\epsilon,M-2}+ C^{\epsilon,M-1})\\  \tau/2 B_3 (C^{\epsilon,M-1}+ C^{\epsilon,M})\\ G^M \\0\\0\\ \vdots\\0\\0\end{array}\right)\,,
\end{equation*}
where the block matrices $BC, BG$ and $Bf$ are given by
\[BC := \left(\begin{array}{cccccc} I&  \mathbf{0}& \cdots& \mathbf{0}&\mathbf{0}& \mathbf{0}\\
A_2&A_1&\cdots&\mathbf{0}&\mathbf{0}& \mathbf{0}\\
\mathbf{0}&\ddots & \ddots & \mathbf{0} &\mathbf{0}& \mathbf{0}\\
\mathbf{0}&\mathbf{0} & \ddots & A_2 & A_1& \mathbf{0}\\
\mathbf{0}&\mathbf{0}&\cdots& \mathbf{0}& A_2 & A_1\\
\frac{\tau}{2} B_3&\frac{\tau}{2} B_3&\cdots&\mathbf{0}&\mathbf{0}&\mathbf{0}\\
\mathbf{0}&\ddots & \ddots &\mathbf{0} & \mathbf{0}&\mathbf{0}\\
\mathbf{0}&\mathbf{0}&\cdots& \frac{\tau}{2} B_3 & \frac{\tau}{2} B_3&\mathbf{0}\\
\mathbf{0}&\mathbf{0}&\cdots&\mathbf{0}& \frac{\tau}{2} B_3 & \frac{\tau}{2} B_3\\
\mathbf{0}&\mathbf{0}&\cdots&\mathbf{0}&\mathbf{0}&\mathbf{0}\\
\mathbf{0}&\mathbf{0}&\cdots&\mathbf{0}&\mathbf{0}&\mathbf{0}\\
\mathbf{0}&\mathbf{0}&\cdots&\mathbf{0}&\mathbf{0}&\mathbf{0}\\
\mathbf{0}&\mathbf{0}&\cdots&\mathbf{0}&\mathbf{0}&\mathbf{0}\\
\mathbf{0}&\mathbf{0}&\cdots&\mathbf{0}&\mathbf{0}&\mathbf{0}\\
\mathbf{0}&\mathbf{0}&\cdots&\mathbf{0}&\mathbf{0}&\mathbf{0}\end{array}\,,\right)\]
\[BG := \left(\begin{array}{ccccccc} \mathbf{0}&\mathbf{0}& \mathbf{0}& \cdots& \mathbf{0}&\mathbf{0}&\mathbf{0}\\
\mathbf{0}& \mathbf{0}& \mathbf{0}& \cdots&\mathbf{0}&\mathbf{0}&\mathbf{0}\\
\mathbf{0}& \mathbf{0}& \mathbf{0}&\cdots&\mathbf{0}&\mathbf{0}&\mathbf{0}\\
\mathbf{0}& \mathbf{0}& \mathbf{0}&\cdots&\mathbf{0}&\mathbf{0}&\mathbf{0}\\
\mathbf{0}& \mathbf{0}& \mathbf{0}&\cdots&\mathbf{0}&\mathbf{0}&\mathbf{0}\\
B_1&B_2&\mathbf{0}& \cdots&\mathbf{0}&\mathbf{0}&\mathbf{0}\\
\mathbf{0}&\ddots & \ddots &\cdots&\mathbf{0}& \mathbf{0}&\mathbf{0}\\
\mathbf{0}& \mathbf{0}& \mathbf{0}&\cdots&B_1&B_2&\mathbf{0}\\
\mathbf{0}& \mathbf{0}& \mathbf{0}&\cdots&\mathbf{0}&B_1&B_2\\
\mathbf{0}& \mathbf{0}& \mathbf{0}&\cdots&\mathbf{0}&\mathbf{0}&I\\
-D^{00}& -D^{01}& \mathbf{0}&\cdots&\mathbf{0}&\mathbf{0}&\mathbf{0}\\
-D^{10}& -D^{11}& -D^{12}&\cdots&\mathbf{0}&\mathbf{0}&\mathbf{0}\\
\mathbf{0}&\ddots& \ddots& \ddots& \mathbf{0}&\mathbf{0}&\mathbf{0}\\
\mathbf{0}& \mathbf{0}& \mathbf{0}&\cdots&-D^{M-1,M-2}&-D^{M-1,M-1}&-D^{M-1,M}\\
\mathbf{0}& \mathbf{0}& \mathbf{0}&\cdots&\mathbf{0}&-D^{M,M-1}&-D^{MM}\end{array}\right)\]
\[Bf := \left(\begin{array}{ccccccc} \mathbf{0}&\mathbf{0}& \mathbf{0}& \cdots& \mathbf{0}&\mathbf{0}&\mathbf{0}\\
-\frac{\tau}{2} B&-\frac{\tau}{2} B&\mathbf{0}&\cdots& \mathbf{0}&\mathbf{0}&\mathbf{0}\\
\mathbf{0}& \ddots & \ddots  &\cdots&\mathbf{0}&\mathbf{0}&\mathbf{0}\\
\mathbf{0}& \mathbf{0} & \ddots  &\cdots&-\frac{\tau}{2} B&-\frac{\tau}{2} B&\mathbf{0}\\
\mathbf{0}&\mathbf{0}&\mathbf{0}&\cdots&\mathbf{0}&-\frac{\tau}{2} B&-\frac{\tau}{2} B\\
\mathbf{0}&\mathbf{0}&\mathbf{0}&\cdots&\mathbf{0}&\mathbf{0}&\mathbf{0}\\
\mathbf{0}&\mathbf{0}&\mathbf{0}&\cdots& \mathbf{0}&\mathbf{0}&\mathbf{0}\\
\mathbf{0}&\mathbf{0}&\mathbf{0}&\cdots&\mathbf{0}&\mathbf{0}&\mathbf{0}\\
\mathbf{0}&\mathbf{0}& \mathbf{0}&\cdots&\mathbf{0}&\mathbf{0}&\mathbf{0}\\
\mathbf{0}&\mathbf{0}& \mathbf{0}&\cdots&\mathbf{0}&\mathbf{0}&\mathbf{0}\\
W^{00}& W^{01}& \mathbf{0}&\cdots&\mathbf{0}&\mathbf{0}&\mathbf{0}\\
W^{10}& W^{11}& W^{12}&\cdots&\mathbf{0}&\mathbf{0}&\mathbf{0}\\
\mathbf{0}&\ddots& \ddots& \ddots&\mathbf{0}&\mathbf{0}&\mathbf{0}\\
\mathbf{0}&\mathbf{0}& \mathbf{0}& \cdots&W^{M-1,M-2}&W^{M-1,M-1}&W^{M-1,M}\\
\mathbf{0}&\mathbf{0}&\mathbf{0}&\cdots&\mathbf{0} &W^{M,M-1}&W^{MM}\end{array}\right)\,.\]

\begin{thebibliography}{}
 \bibitem{ACM99} Aitbayev, R., Cai,  X.-C., Paraschivoiu, M.: Parallel two-level methods for three-dimensional transonic compressible flow simulations on unstructured meshes. Proceedings of Parallel CFD'99 (1999)
 \bibitem{AB05} Akcelik, V., Biros, G., Draganescu, A., Ghattas, O., Hill, J., Waanders, B.: Dynamic data-driven inversion for terascale simulations: Real-time identification of airborne contaminants. Proceedings of Supercomputing, Seattle, WA (2005)
\bibitem{AB03} Akcelik, V., Biros, G., Ghattas, O., Long, K.~R., Waanders, B.: A variational finite element method for source inversion for convective-diffusive transport. Finite Elem. Anal. Des. \textbf{39}, 683-705 (2003)
\bibitem{AB01} Atmadja, J., Bagtzoglou, A.~C.: State of the art report on mathematical methods for groundwater pollution source identification. Environ. Forensics \textbf{2}, 205-214 (2001)
\bibitem{BBMTZ02} Baflico, L., Bernard, S., Maday, Y., Turinici, G., Zerah, G.: Parallel-in-time molecular-dynamics simulations. Phys. Rev. E \textbf{66},  2-5 (2002)
\bibitem{BBEG10} Balay, S., Buschelman, K., Eijkhout, V., Gropp, W.~D., Kaushik, D., Knepley, M.~G., McInnes, L.~C., Smith, B.~F., Zhang, H.: PETSc Users Manual. Technical Report, Argonne National Laboratory (2014)
\bibitem{B96} Battermann, A.: Preconditioners for Karush-Kuhn-Tucker Systems Arising in Optimal Control. Master Thesis, Virginia Polytechnic Institute and State University, Blacksburg, Virginia (1996)
\bibitem{BG99}	Biros, G., Ghattas, O.: Parallel preconditioners for KKT systems arising in optimal control of viscous incompressible flows.  Proceedings of Parallel CFD'99, Williamsburg, Virginia, USA (1999)
\bibitem{CLZ09} Cai, X.-C.,  Liu, S., Zou, J.: Parallel overlapping domain decomposition methods for coupled inverse elliptic problems. Comm. App. Math. Com. Sc. \textbf{4}, 1-26 (2009)
\bibitem{CS99} Cai, X.-C., Sarkis, M.: A restricted additive Schwarz preconditioner for general sparse linear systems. SIAM J. Sci. Comput. \textbf{21},  792-797 (1999).
\bibitem{CC12} Chen, R.~L., Cai, X.-C.: Parallel one-shot Lagrange-Newton-Krylov-Schwarz algorithms for shape optimization of steady incompressible flows. SIAM J. Sci. Comput. \textbf{34}, 584-605 (2012)
\bibitem{DZZ13} Deng, X.~M., Zhao, Y.~B., Zou, J.: On linear finite elements for simultaneously recovering
    source location and intensity.  Int. J. Numer. Anal. Mod. \textbf{10}, 588-602 (2013)
\bibitem{EHN98}  Engl, H.~W., Hanke, M., Neubauer, A.: Regularization of Inverse Problems. Kluwer Academic Publishers, Netherland (1998)
\bibitem{FC03} Farhat, C., Chandesris, M.: Time-decomposed parallel time-integrators: theory and feasibility studies for fluid, structure, and fluid-structure applications.  Int. J. Numer. Meth. Eng. \textbf{58}, 1397-1434 (2003)
\bibitem{GH08}  Gander,  M.~J., Hairer, E.: Nonlinear convergence analysis for the parareal algorithm. Proceedings of the 17th International Conference on Domain Decomposition Methods \textbf{60}, 45-56 (2008)
\bibitem{GP08} Gander, M.~J., Petcu, M.: Analysis  of  a  Krylov  subspace  enhanced  parareal algorithm for linear problems.  Paris- Sud Working Group on Modeling and Scientific Computing 2007- 2008 (E. Cances et al., eds.),  ESAIM Proc. EDP Sci., LesUlis \textbf{25}, 114-129 (2008)
\bibitem{GV07} Gander, M.~J., Vandewalle, S.: Analysis of the parareal time-parallel time-integration method. SIAM J. Sci. Comput. \textbf{29}, 556-578 (2007)
\bibitem{GER83}	Gorelick, S., Evans, B., Remson, I.: Identifying sources of groundwater pollution: an optimization approach. Water Resour. Res. \textbf{19}, 779-790 (1983)
\bibitem{H09}   Hamdi, A.: The recovery of a time-dependent point source in a linear transport equation: application to surface water pollution. Inverse Probl., \textbf{24}, 1-18 (2009)
\bibitem{KGM11} Karalashvili, M., Gro$\beta$, S., Marquardt, W., Mhamdi, A., Reusken, A.: Identification of transport coefficient models in convection-diffusion equations. SIAM J. Sci. Comput. \textbf{33}, 303-327 (2011)
\bibitem{KZ98}  Keung, Y.~L., Zou, J.: Numerical identifications of parameters in parabolic systems. Inverse Probl. \textbf{14}, 83-100 (1998)
\bibitem{KT51}  Kuhn, H.~W., Tucker, A.~W.: Nonlinear programming. Proceedings of 2nd Berkeley Symposium, Berkeley: University of California Press, 481-492 (1951)
\bibitem{LMT01} Lions, J.-L., Maday, Y., Turinici, G.: A ¡°parareal¡± in time discretization of PDE's. ComptesRendus de l'Academie des Sciences Series I Mathematics \textbf{332}, 661-668 (2001)
\bibitem{LZ07}	Liu, X., Zhai, Z.: Inverse modeling methods for indoor airborne pollutant tracking literature review and fundamentals. Indoor Air \textbf{17},  419-438 (2007)
\bibitem{MD09} Zhang, J., Delichatsios, M.~A.: Determination of the convective heat transfer coefficient in three-dimensional inverse heat conduction problems.  Fire Safety J. \textbf{44}, 681-690 (2009)
\bibitem{MT05} Maday, Y., Turinici G.: The parareal in time iterative solver: a further direction to parallel implementation. Domain Decomposition Methods in Science and Engineering, Springer LNCSE \textbf{40}, 441-448 (2005)
\bibitem{NKM09} Nilssen, T.~K., Karlsen, K.~H., Mannseth, T., Tai, X.-C.: Identification of diffusion parameters in a nonlinear convection-diffusion equation using the augmented Lagrangian method. Computat. Geosci.  \textbf{13}, 317-329 (2009)
\bibitem{PBC06} Prudencio, E., Byrd, R., Cai, X.-C.: Parallel full space SQP Lagrange-Newton-Krylov-Schwarz algorithms for PDE-constrained optimization problems. SIAM J. Sci. Comput. \textbf{27}, 1305-1328 (2006)
\bibitem{RR05}  Revelli, R., Ridolfi, L.: Nonlinear convection-dispersion  models with  a  localized  pollutant  source II--a  class  of  inverse  problems. Math. Comput. Model. \textbf{42}, 601-612 (2005)
\bibitem{S93}  Saad, Y.: A flexible inner-outer preconditioned GMRES algorithm. SIAM J. Sci. Comput.  \textbf{14}, 461-469 (1993)
\bibitem{SV07} Samarskii, A.~A., Vabishchevich, P.~N.: Numerical Methods for Solving Inverse Problems of Mathematical Physics. Walter de Gruyter, Berlin  (2007).
\bibitem{SK94}  Skaggs, T., Kabala, Z.: Recovering the release history of a groundwater contaminant.  Water Resour. Res. \textbf{30},  71-80 (1994)
\bibitem{SK95}  Skaggs, T., Kabala, Z.: Recovering the history of a groundwater contaminant plume: method of quasi-reversibility.  Water Resour. Res. \textbf{31},  2669-2673 (1995)
\bibitem{SBG04}  Smith, B., Bj{\o}rstad, P., Gropp, W.: Domain Decomposition: Parallel Multilevel Methods for Elliptic Partial Differential Equations. Cambridge University Press (2004)
\bibitem{SK97}	Snodgrass, M.~F., Kitanidis, P.~K.: A geostatistical approach to contaminant source identification. Water Resour. Res. \textbf{33}, 537-546 (1997)
\bibitem{WY10}  Wong, J., Yuan, P.: A FE-based algorithm for the inverse natural convection problem.  Int. J. Numer. Meth. Fl., \textbf{68}, 48-82 (2012)
\bibitem{W03}   Woodbury, K.~A.: Inverse Engineering Handbook. CRC Press (2003)
\bibitem{YPC12} Yang, H., Prudencio, E., Cai, X.-C.: Fully implicit Lagrange-Newton-Krylov-Schwarz algorithms for boundary control of unsteady incompressible flows. Int. J. Numer. Meth. Eng. \textbf{91}, 644-665 (2012)
\bibitem{YSL08} Yang, X.-H., She, D.-X., Li, J.-Q.: Numerical approach to the inverse convection-diffusion problem. 2007 International Symposium on Nonlinear Dynamics (2007 ISND), Journal of Physics: Conference Series \textbf{96},  012156 (2008)
\end{thebibliography}
\end{document}